\documentclass[amsfonts,12pt]{article}
\def \ed{\end{document}}

\usepackage{latexsym,epsfig,graphpap,graphics,amsmath}
\usepackage{amssymb,amsmath,latexsym,graphicx}
\numberwithin{equation}{section}
\def \n1{\newpage}
\def \L1{\frak L}

\def \p{\partial}

\def \bqn{\begin{equation}}
\def \9{\end{equation}}
\def \3{\begin{eqnarray*}}
\def \4{\end{eqnarray*}}
\def \1{\begin{eqnarray}}
\def \2{\end{eqnarray}}

\def \big1{\bigcap}

\def \fr{\frac}

 \newcounter{corollary}

\newcounter{proposition}
\newcounter{definition}

\def \brem{\begin{remark}}
\def \erem{\end{remark}}
\def \bth{\begin{theorem}}
\def \eth{\end{theorem}}
\def \bpr{\begin{proposition}}
\def \epr{\end{proposition}}
\def \bprf{\begin{proof}}
\def \eprf{\end{proof}}
\def \bex{\begin{example}}
\def \eex{\end{example}}
\def \bprf{\begin{proof}}
\def \eprf{\end{proof}}
\def \blem{\begin{lemma}}
\def \elem{\end{lemma}}
\newcounter{theorem}

\def \big{\bigcap}

\def \bl{\begin{lemma}}
\def \bcor{\begin{corollary}}
\def \ecor{\end{corollary}}
\def \el{\end{lemma}}
\def \beq*{\begin{eqnarray*}}
\def \eeq*{\end{eqnarray*}}
\def \6{\vspace*{7mm}}

\def \s1{\sqrt}

\def \mb{\mbox}

\def \bt{\begin{tabular}}
\def \et{\end{tabular}}

\def \l{\left}

\def \r{\right}
\def \hs1{\hspace*{3mm}}
\def \q2{\hspace*{9mm}}

\def \un1{\underline}
\def \mb{\mbox}
\def \vs1{\vspace*{4mm}}

\def \ba{\begin{array}}
\def \ea{\end{array}}

\newcommand{\ec}{\end{center}}
\newcommand{\bc}{\begin{center}}
\newcommand{\be}{\begin{equation}}
\newcommand{\ds}{\displaystyle}

\newcommand{\ee}{\end{equation}}
\newcommand{\bn}{\begin{enumerate}}
\newcommand{\en}{\end{enumerate}}
\newcommand{\bi}{\begin{itemize}}
\newcommand{\ei}{\end{itemize}}
\usepackage{graphicx}

\newtheorem{theorem}{Theorem}

\newtheorem{corollary}{Corollary}

\newtheorem{example}{Example}

\newtheorem{lemma}{Lemma}

\newtheorem{proposition}{Proposition}
\newtheorem{remark}{Remark}

\newenvironment{proof}[1][Proof]{\textbf{#1.} }{\
\rule{0.5em}{0.5em}}
\graphicspath{{converted_graphics/}}
\begin{document}
\title{Analytic solutions of initial-boundary-value problems of transient conduction using symmetries}
\author{H. Azad$^{1}$, M. T. Mustafa$^{1}$  and A. F. M. Arif$^{2}$}
\date{}
\maketitle \vspace{-.5in}
\begin{center}
$^{1}$Department of Mathematics and Statistics, \\ King Fahd University of Petroleum and Minerals, \\
Dhahran 31261, Saudi Arabia\\
 \texttt{hassanaz@kfupm.edu.sa}, \texttt{tmustafa@kfupm.edu.sa}\\
$^{2}$Department of Mechanical Engineering, \\ King Fahd University of Petroleum and Minerals, \\
Dhahran 31261, Saudi Arabia\\
\texttt{afmarif@kfupm.edu.sa}\\
\end{center}

\baselineskip=18pt
\begin{abstract}
\noindent Lie symmetry method is applied to find analytic
solutions of initial-boundary-value problems of transient
conduction in semi-infinite solid with constant surface
temperature or constant heat flux condition. The solutions are
obtained in a manner highlighting the systematic procedure of
extending the symmetry method for a PDE to investigate BVPs of the
PDE. A comparative analysis of numerical and closed form solutions
is carried out for a physical problem of heat conduction in a
semi-infinite solid bar made of AISI 304 stainless steel.
\end{abstract}

%\vs1
%
%\noindent 2000 Mathematics subject classification:
%                        \\[2mm]
%Key words:

\renewcommand{\theequation}{\thesection.\arabic{equation}}
\section{Introduction}

Lie symmetry method is a powerful general technique for analyzing
non-linear PDEs and can be efficiently employed to study problems
having implicit or explicit symmetries. Thus it provides most
widely applicable technique to find the closed form solutions of
differential equations and contains, as particular case cf.
\cite{pz04}, many efficient methods for solving differential
equations like separation of variables, traveling wave solutions,
self-similar solutions and exponential self-similar solutions.
Since the modern treatment of the classical Lie symmetry theory by
Ovsiannikov \cite{ovsiannikov}, the theory of symmetries of
differential equations has been studied intensely and has
substantially grown. A large amount of literature about the
classical Lie symmetry theory, its applications and its extensions
is available, e.g.
\cite{Ames65-72,baumann,bluman-cole,bluman-kumei,euler,hansen,hydon,Ib-crc1,Ib-crc2,Ib-crc3,Ib-book,miller77,olver,ovsiannikov,stephani}.

Most of the engineering and physical problems require the PDE to
be solved subject to suitable initial or boundary conditions.
Although there have been some notable contributions in the
applications of symmetry method to boundary value problems (BVPs)
c.f. \cite{Ames65-72, bluman-bvp, cantwell, dresner, ovsiannikov},
in general the Lie symmetry method has not been utilized in a
great deal in obtaining solutions of BVPs of physical
significance. One reason could be the natural restrictions imposed
because of the requirement of invariance of initial-boundary
conditions in addition to requirement of invariance of PDE under
the symmetry. This work is concerned with application of Lie
symmetry method to obtain analytic solution of two standard
initial-boundary-value problems of heat conduction, in a
systematic manner that highlights the systematic procedure of
extending the symmetry method for a PDE to investigate BVPs of the
PDE.

A description of the initial-boundary value problems, under study,
is provided in Section~2. Section~3 is divided in two parts,
giving an introduction to the solution-method followed by its
application to obtain exact solutions of the initial-boundary
value problems described in Section~2. The comparative analysis of
numerical and analytic solutions for the heat conduction in solid
bar made of AISI 304 stainless steel is presented in Section~4.

%%%%%%%%%%%%%%%%%%%%%%%%%%%%%%%%%%%%%%
%
%
%%%%%%%%%%%%%%%%%%%%%%%%%%%%%%%%%%%%%%
\section{Description of the initial-boundary value problems}

We consider test problems related to transfer of heat by
conduction. The analysis of such problems is required in many
physical engineering problems, for example, the cooling of
electronic equipment, the design of thermal-fluid systems, and the
material and manufacturing processes. In practice, a major
objective of the solution of such problems is to determine the
temperature field in a medium as a result of either thermal
condition applied to the boundary of the medium or heat generation
within the medium. Once the exact temperature distribution is
known, the heat flux at any point in the system, including the
boundaries, can be computed form the Fourier's law.

The problems studied here are for the transient conduction in
semi-infinite solid either with constant surface temperature or
with constant surface heat flux conditions. Precisely, we
investigate the following initial-boundary value problems.

\begin{description}
\item [IBVP-1]\hfill \\ {\bf Transient conduction in a
semi-infinite solid with constant surface temperature}
\begin{equation}\label{1}
\frac{\partial T}{\partial t}=\alpha \frac{\partial
^{2}T}{\partial x^{2}},
\end{equation}
with initial and boundary conditions
\begin{equation}\label{2}
T|_{t=0} = T_i,\qquad T|_{x=0} = T_s, \qquad T|_{x\to \infty} =
T_i
\end{equation}
\item [IBVP-2] \hfill \\ {\bf Transient conduction in a
semi-infinite solid with constant surface heat flux}
\end{description}
\begin{equation}\label{3}
\frac{\partial T}{\partial t}=\alpha \frac{\partial
^{2}T}{\partial x^{2}},
\end{equation}
with initial and boundary conditions
\begin{equation}\label{4}
T|_{t=0} = 0,\qquad -k \frac{\partial T}{\partial x}|_{x=0}=
{q_0}^{\prime\prime}, \qquad T|_{x\to \infty} = 0
\end{equation}
%
%%%%%%%%%%%%%%%%%%%%%%%%%%%%%%%%%%%%%%
%
%
%%%%%%%%%%%%%%%%%%%%%%%%%%%%%%%%%%%%%%
\section{Lie symmetry solution of the initial-boundary value problems}

The general procedure of applying Lie symmetry method to study
IBVPs of PDEs requires \cite{bluman-anco} determination of a one
parameter Lie group of transformations that leaves the problem
invariant, and the utilization of these transformations either to
construct the invariant solution or to obtain similarity
reductions. It can be explained by the following steps:
\begin{itemize}
\item {\it Lie symmetries of the PDE} \\ Determining the symmetry
algebra of the governing PDE.
\item {\it Invariance of the boundaries} \\ Taking the most
general symmetry operator $X$ obtained in step~1 and finding the
conditions under which it leaves the boundaries invariant.
\item {\it Invariance of boundary conditions on the boundary} \\
Finding the restrictions on $X$ that are imposed due to invariance
of boundary conditions on the boundary.
\end{itemize}
These steps will determine the symmetry operator that leaves the
IBVP invariant.
\begin{itemize}
\item {\it Construction of similarity solution or reductions} \\
Utilizing the similarity variables of the symmetry operator of the
IBVP to find the similarity reductions and similarity solution of
IBVP.
\end{itemize}
Further details about these steps are illustrated in the
subsections below. It should be noted that the more general the
symmetry operator (leaving the IBVP invariant) found the more
likely it is to lead to the solution of the problem \cite{Ib-pc}.

\subsection{Lie symmetries of PDE (\ref{1})}
The method of obtaining the classical Lie symmetries of a PDE is
standard which is described in detail in many books, e.g.
\cite{bluman-kumei, hydon, Ib-book, olver, stephani}. For
instance, to obtain Lie symmetries of an equation like (\ref{1}),
one considers the 1-parameter Lie group of infinitesimal
transformations in $(x,t,T)$ given by
\[
\renewcommand{\arraystretch}{2}
\ba{l}
x^* = x+\epsilon \xi (x,t,T)+ O(\epsilon^2) \\
t^* = t+\epsilon \tau (x,t,T)+ O(\epsilon^2) \\
T^*=  T+\epsilon \phi (x,t,T)+ O(\epsilon^2) \ea \] where
$\epsilon$ is the group parameter, hence the corresponding
generator of the Lie algebra is of the form
\[
X = \xi(x,t,T) \fr{\p}{\p x} + \tau (x,t,T) \fr{\p}{\p t}
     + \phi (x,t,T)\fr{\p}{\p T}. \]
If $X^{[2]}$ denotes the second prolongation of $X$ then using the
invariance condition
$$
X^{[2]}\left( T_{t}-\alpha T_{xx}  \right)|_{T_{t}=\alpha
T_{xx}}=0
$$
yields an overdetermined system of linear PDEs in $\xi$, $\tau$ and
$\phi$ called determining equations. The general solution of this
system determines the generator of the symmetry algebra. For most
type of problems, Computer Algebra Systems like Macsyma, Maple,
Mathematica, MuPAD and Reduce can be used either to perform above
steps separately or to find the symmetry algebra directly. A
detailed survey of software that can be used for symmetry analysis
is given by Hereman in \cite[pages 367-413]{Ib-crc3}.

The symmetry algebra of PDE~(\ref{1}) is well known (in fact it was
found by Lie) and is spanned by the vector fields
$$
X_1 =\fr{\p}{\p t},\quad X_2 =\fr{\p}{\p x},\quad X_3
=2t\fr{\p}{\p t} + x\fr{\p}{\p x},\quad X_4 =2t\fr{\p}{\p x} -
\frac{1}{\alpha} xT\fr{\p}{\p T},
$$
$$
X_5 =4t^2 \fr{\p}{\p t} +4xt\fr{\p}{\p x} -\frac{1}{\alpha}
(x^2+2\alpha t)T\fr{\p}{\p T},\quad X_6 =T\fr{\p}{\p T}, \quad
X_{\infty} =f(t,x)\fr{\p}{\p T}.
$$
The corresponding optimal system of 1-dimensional subalgebras is
given by the following vector fields, \cite{olver}.
$$
X_2,\quad X_6,\quad X_1+cX_6,\quad X_1 +X_4,\quad X_1 - X_4,\quad
X_1+X_5+cX_6,\quad X_3 +cX_6
$$
and representatives of the form $X_\infty$. This means any invariant
solution of PDE~(\ref{1}) can be found via a suitable transformation
of the invariant solutions obtained from the symmetry operators in
the optimal system. Instead, for our purpose we directly find the
symmetry operator that preserves the boundary and boundary
conditions and hence directly leads to the solutions of
initial-boundary-value problems (see details below).
%%%%%%%%%%%%%%%%%%%%%%%%%%%%%%%
%%
%%%%%%%%%%%%%%%%%%%%%%%%%%%%%%%
\subsection{Solution to IBVP-1}
We consider the general symmetry operator
\[
X = k_1 X_1 + k_2 X_2 + k_3 X_3 + k_4 X_4 + k_5 X_5 + k_6 X_6 \] of
PDE~(\ref{1}) and search for the operator that preserves the
boundary and the boundary conditions~(\ref{2}).

The invariance of the boundaries $x=0, t = 0$ or equivalently \3
[X(x-0)]_{x=0} & = & 0, \\
\mb{}[X(t-0)]_{t=0} & = & 0 \4 implies
\bqn\label{cond1-bvp1} k_1 = k_2 = k_4 = 0. \9
Hence, $X$ must be $$X = k_3 X_3  + k_5 X_5 + k_6 X_6 .$$ In
addition to the restrictions imposed by
Equation~(\ref{cond1-bvp1}), the invariance of initial and
boundary conditions i.e. \3
 [X(T-T_i)]_{t=0} & = & 0, \quad {\rm on}\ T=T_i \\
 \  [X(T-T_s)]_{x=0} & = & 0,  \quad {\rm on}\ T=T_s \4
implies we must have \bqn\label{cond2-bvp1} k_5 = 0 = k_6. \9
Hence the IBVP (\ref{1}), (\ref{2}) is invariant under the
symmetry
\[
X = 2t \fr{\p}{\p t} + x \fr{\p }{\p x},\] where we have chosen $k_3
=1$.

The invariant solution of the problem is constructed by utilizing
the transformations through similarity variables for $X$. Solving
the characteristic system for $XI = 0$ gives $\ds I_1 =
\fr{x^2}{t}$ and $I_2 = T$ as the differential invariants of $\ds
X = 2t\fr{\p}{\p t} + x \fr{\p }{\p x}$. Hence, the similarity
variables for $X$ are \bqn \label{bvp1-sim-var} \xi(x,t) =
\fr{x^2}{t} \mb{ and } V(\xi) = T. \9 Substitution of similarity
variables in Equation~(\ref{1}) implies that the corresponding
similarity solution of PDE~(\ref{1}) is of the form $T = V(\xi)$
where $V(\xi)$ satisfies the ODE \bqn \label{bvp1-ode-xi} 4\xi
\fr{d^2V}{d\xi^2} + \l(2 + \fr{\xi}{\alpha}\r)\fr{dV}{d\xi} = 0.
\9
The above equation can be integrated, using the substitution
\[
W = \fr{dV}{d\xi},\] to obtain \bqn \label{bvp1-sol-xi} V(\xi) =
c_1 \int \fr{e^{-\xi/4\alpha}}{\s1{\xi}} d\xi + c_2. \9 Making the
change of variable \bqn \label{bvp1-c-var} y^2 = \fr{\xi}{4\alpha}
\9 in the above solution yields \3
V & = & 4c_1 \s1{\alpha} \int e^{-y^2}dy + c_2 \\
& = & 2c_1 \s1{\alpha} \s1{\pi}\cdot \mathbf{erf} (y) + c_2, \4
where $\mathbf{erf}$ denotes the error function. Hence, from
Equations~(\ref{bvp1-sim-var}), (\ref{bvp1-sol-xi}) and
(\ref{bvp1-c-var}), the exact solution of PDE (\ref{1}) that is
invariant under $\ds X = 2t \fr{\p}{\p t} + x\fr{\p}{\p x}$ is
\bqn T(x,t) = 2c_1 \s1{\alpha} \s1{\pi} \cdot \mathbf{erf} \l(
\fr{x/\s1{t}}{2\s1{\alpha}}\r)+c_2. \9 Imposing the initial and
boundary conditions determines
\[
c_1 = \fr{T_i - T_s}{2\s1 {\alpha\pi} } \mb{ and } c_2 = T_s,\]
giving the solution
\[ T(x,t) = (T_i - T_s) \mathbf{erf} \l( \fr{x}{2\s1{\alpha t}}\r) + T_s\]
of the IBVP (\ref{1}), (\ref{2}).

%%%%%%%%%%%%%%%%%%%%%%%%%%%%%%%%%%%%%%%%%%%%%%%%%%%%
%%
%%%%%%%%%%%%%%%%%%%%%%%%%%%%%%%%%%%%%%%%%%%%%%%%%%%%
\subsection{Solution to IBVP-2}
Taking the general symmetry operator
\[
X = k_1 X_1 + k_2 X_2 + k_3 X_3 + k_4 X_4 + k_5 X_5 + k_6 X_6 \]
of PDE~(\ref{3}), we determine the operator that leaves the
boundary and the boundary conditions~(\ref{4}) invariant.

The invariance of the boundaries $x=0, t = 0$ or equivalently \3
[X(x-0)]_{x=0} & = & 0, \\
\mb{}[X(t-0)]_{t=0} & = & 0 \4 yields
\bqn\label{cond1-bvp2} k_1 = k_2 = k_4 = 0, \9 i.e. $X$ must be
\[
X = k_3 X_3  + k_5 X_5 + k_6 X_6. \]
Since
\begin{equation}
[X(T-0)]_{t=0} = 0
\end{equation}
implies
\begin{equation}
\l( \fr{-k_5}{\alpha} x^2 + k_6 \r) T = 0,
\end{equation}
the invariance of the condition $T|_{t=0} = 0$ at $T=0$ does not
impose any restrictions on $X$.

The invariance of the boundary condition
$$
\l[-k \frac{\partial T}{\partial x}= {q_0}^{\prime\prime}\r]_{x=0}
$$
requires, see \cite[Chapter 4]{bluman-anco} for details,
$$
\l[ X^{[1]}\l( k \frac{\partial T}{\partial x}+
{q_0}^{\prime\prime} \r) \r]_{x=0}=0 \qquad {\rm on}\ -k
\frac{\partial T}{\partial x}= {q_0}^{\prime\prime}
$$
where $X^{[1]}$ denotes the first prolongation of $X$. It follows
that
$$
\l[ \l( -6k_5 t + k_6-k_3 \r) \frac{\partial T}{\partial x}
\r]_{x=0} =0
$$
on
$$
k \frac{\partial T}{\partial x}+ {q_0}^{\prime\prime}=0,
$$
hence we must have
$$
k_5=0 \mb{  and  } k_3 = k_6.
$$
Choosing $k_3 = k_6 = 1$ provides the symmetry
$$
X= X_3 + X_6 = x \fr{\p }{\p x} + 2 t \fr{\p}{\p t} + T \fr{\p}{\p
T}
$$
that leaves the IBVP (\ref{3}), (\ref{4}) invariant.

To find the similarity transformations that will lead to the
solution, the characteristics system for $XI = 0$ is solved. This
provides the similarity variables \bqn \label{bvp2-sim-var}
\xi(x,t) = \fr{x^2}{t} \mb{ and } V(\xi) = \frac{T}{x}. \9
Substituting similarity variables in PDE~(\ref{3}) implies that
the corresponding similarity solution is of the form \bqn
\label{bvp2-sol-form} T = x V(\xi) \9 where $V(\xi)$ satisfies the
ODE
\bqn \label{bvp2-ode-xi} 4 \fr{d^2V}{d\xi^2} + \l(\frac{6}{\xi} +
\fr{1}{\alpha}\r)\fr{dV}{d\xi} = 0. \9 Following a procedure
similar to the solution of ODE~(\ref{bvp1-ode-xi}) and using
Equations~(\ref{bvp2-sim-var}), (\ref{bvp2-sol-form}) we obtain
\bqn T(x,t) =-2c_1 \s1{t}e^{\frac{-x^2}{4\alpha t}}-c_1 x
\s1{\frac{\pi}{\alpha}}\cdot
 \mathbf{erf} \l(\fr{x/\s1{t}}{2\s1{\alpha}}\r) + c_2 x . \9
Imposing the boundary conditions determines
\[
c_1 = -\frac{{q_0}^{\prime\prime}}{\s1 {k}}
\s1{\frac{\alpha}{\pi}} \mb{ and } c_2 =
-\frac{{q_0}^{\prime\prime}}{k} .\]
giving the solution
$$
T(x,t)=2 \frac{{q_0}^{\prime\prime}}{k} \s1{\frac{\alpha t}{\pi}}
e^{\frac{-x^2}{4\alpha t}} + \frac{{q_0}^{\prime\prime}}{k} x
\left\{   \mathbf{erf} \l( \frac{x}{2\s1{\alpha t}}  \r) -1
\right\}
$$
 of the IBVP (\ref{3}), (\ref{4}).

\subsubsection*{Remarks}
\begin{enumerate}
\item Although the above procedure illustrates that the symmetry
algebras admitted by boundary value problems, in general, may not
be rich enough to obtain symmetry solutions or reductions of the
BVP, it is worth noting that sometimes the invariant solutions
obtained from PDE can satisfy boundary conditions that are not
left invariant by the operator which generated the solution.
Examples can easily be constructed from above boundary value
problems. For instance, we can see from above that the invariant
solution
$$
T(x,t)= T_i \mathbf{erf} \l( \fr{x/\s1{t}}{2\s1{\alpha}}\r)
$$
obtained from the generator
\[
X_3 = 2t \fr{\p}{\p t} + x \fr{\p }{\p x},\] satisfies the IBVP
$$
\frac{\partial T}{\partial t}=\alpha \frac{\partial
^{2}T}{\partial x^{2}}
$$
with initial and boundary conditions
$$
T|_{t=0} = T_i,\qquad T|_{x=2\s1{\alpha}} = T_i \mathbf{erf}\l(
\frac{1}{\s1{t}} \r),
$$
but the boundary condition
$$
\l[ T=T_i \mathbf{erf}\l( \frac{1}{\s1{t}} \r)
\r]_{x=2\s1{\alpha}}
$$
is not preserved by the generator \[ X_3 = 2t \fr{\p}{\p t} + x
\fr{\p }{\p x}.\]
\item In case of non-linear problems it is likely that reduced
ODE, obtained through similarity variables, may not be integrable
in terms of known functions. In such cases, the reduced BVP of ODE
either can be solved numerically or a Lie plane analysis, cf.
\cite{dresner}, of the reduced ODE has to be carried out to
understand the solution.
\end{enumerate}

%%%%%%%%%%%%%%%%%%%%%%%%%%%%%%%%%%%%%%
%
%
%%%%%%%%%%%%%%%%%%%%%%%%%%%%%%%%%%%%%%

\section{Comparison of numerical and analytic results for a physical problem}
To compare the exact solution of the two boundary value problems
discussed in the previous sections with numerical results using
finite element method, transient heat conduction in a
semi-infinite solid bar was solved. The solid is made of AISI 304
stainless steel and initially at a temperature $T_i$. The material
properties for AISI 304 steel used in this work are given below:
\begin{eqnarray}
 && {\rm Thermal\  conductivity}\ (k)= 18.2\ W/m\ ^{0}K  \nonumber \\
&&  {\rm Specific\ heat }\ (c)= 0.536\ kJ/kg\ ^{0}K \nonumber \\
&&  {\rm Density }\ (\rho)= 7822\ kg/m^3\  \nonumber \\
&&  {\rm Thermal\ diffusivity }\ (\alpha=k/\rho c)= 4.34E-3\ m^2
/s\ \nonumber
\end{eqnarray}

The spatial temperature distribution at different times was
calculated when the solid is subjected to the following boundary
conditions:
\begin{center}
\begin{tabular}{|l|l|}\hline
{\bf For BVP-1} & {\bf For BVP-2}  \\ \hline
Initial temperature $T_i = 300\ ^{0}K$ & Initial temperature $T_i = 0\ ^{0}K$  \\
Surface temperature $T_s = 900\ ^{0}K$ at $x=0$ & Heat flux $= 5\ kW/m^2$ at $x=0$ \\
\hline
\end{tabular}
\end{center}

In order to numerically solve this problem, a transient thermal
analysis is performed using FEA software ANSYS. The given
structure is modeled using 3-D Conduction Bar Elements (LINK33).
LINK33 is a uniaxial element with the ability to conduct heat
between its nodes. The element has a single degree of freedom,
temperature, at each node point. The conducting bar is applicable
to a steady-state or transient thermal analysis. A refined uniform
mesh is used to model the nonlinear thermal gradient through the
solid. The length of the model is taken as $L=2m$ for BVP-1 and
$L=10m$ for BVP-2 assuming that no significant temperature change
occurs at the interior end point during the time period of
interest. This assumption is validated by the temperature of node
at $x=L$ at the end of the transient analysis.

In the figures given at the end, Figure 1 shows the spatial distribution of temperature using
closed form solution of the BVP-1. The comparison of this solution
with the numerical results using FEA at three different times is
shown in Figure 2. Similarly, the temperature distribution for the
BVP-2 from the exact solution is given in Figure-3. Figure-4
indicates very good agreement between closed form solution and
numerical results for the BVP-2.
\section{conclusion}
This paper presents a systematic approach for applying Lie
symmetry method to investigate boundary value problems of partial
differential equations. Equations for transient conduction in
semi-infinite solid are considered subject to either constant
surface temperature condition or constant heat flux condition.
Those symmetries are determined that leave the partial
differential equation, boundaries and the boundary conditions
invariant. The corresponding similarity transformations are then
utilized to transfer the problem to boundary value problem of
ordinary differential equations which are solved to obtain
explicit analytic solutions in both cases. Finally, the results
are applied to practical problems of heat conduction in a solid
bar made up of AISI 304 stainless steel. For these particular
problems, a comparative analysis of numerical and closed form
solutions is also carried out. The recovery of well-known
solutions of important transient conduction BVPs of engineering
shows that a systematic application of Lie symmetry method to BVPs
is an applicable technique and demonstrates its potential for
investigating BVPs related to practical applications in different
fields of applied sciences and technology.

%%%%%%%%%%%%%%%%%%%%%%%%%%%%%%%%%
%
%
%%%%%%%%%%%%%%%%%%%%%%%%%%%%
\section*{Acknowledgments}
The authors acknowledge the support of King Fahd University of
Petroleum $\&$ Minerals, Dhahran, Saudi Arabia.

%%%%%%%%%%%%%%%%%%%%%%%%
%
%   All figures
%
%%%%%%%%%%%%%%%%%%%%%%%%

\newpage

\begin{figure}[h] % float placement: (h)ere, page (t)op, page (b)ottom, other (p)age
  \centering
  % file name: D:/tahir-d/research/current-work/arif-bvps/first-paper/write-up/to-submit/arxiv/figure1.jpg
  \includegraphics[width=5.67in,height=4.02in,keepaspectratio]{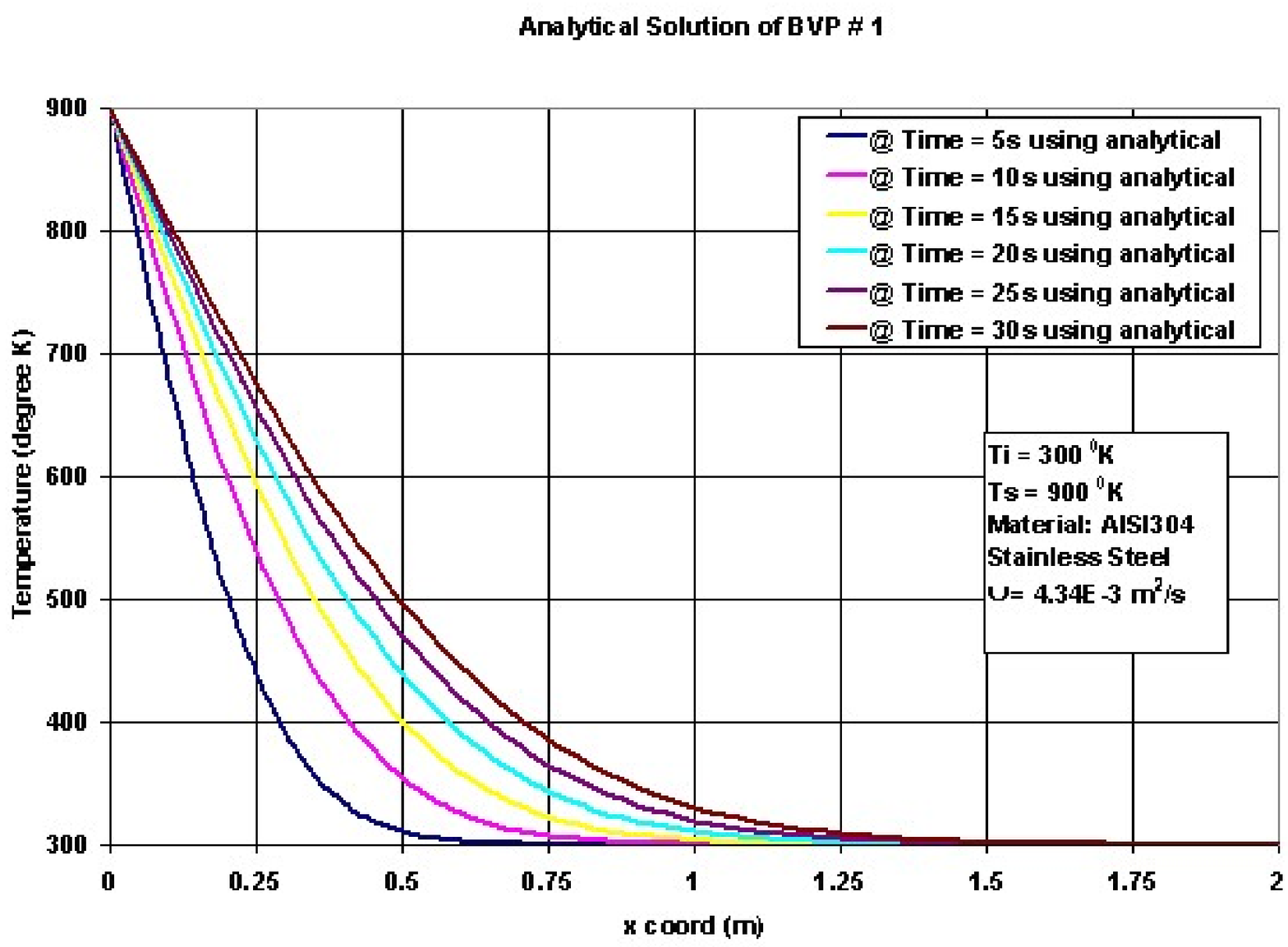}
  \caption{Temperature distribution at different times using analytical solution for BVP-1}
  \label{fig:figure1}
\end{figure}

\begin{figure}[h] % float placement: (h)ere, page (t)op, page (b)ottom, other (p)age
  \centering
  % file name: D:/tahir-d/research/current-work/arif-bvps/first-paper/write-up/to-submit/arxiv/figure2.jpg
  \includegraphics[width=5.67in,height=4.08in,keepaspectratio]{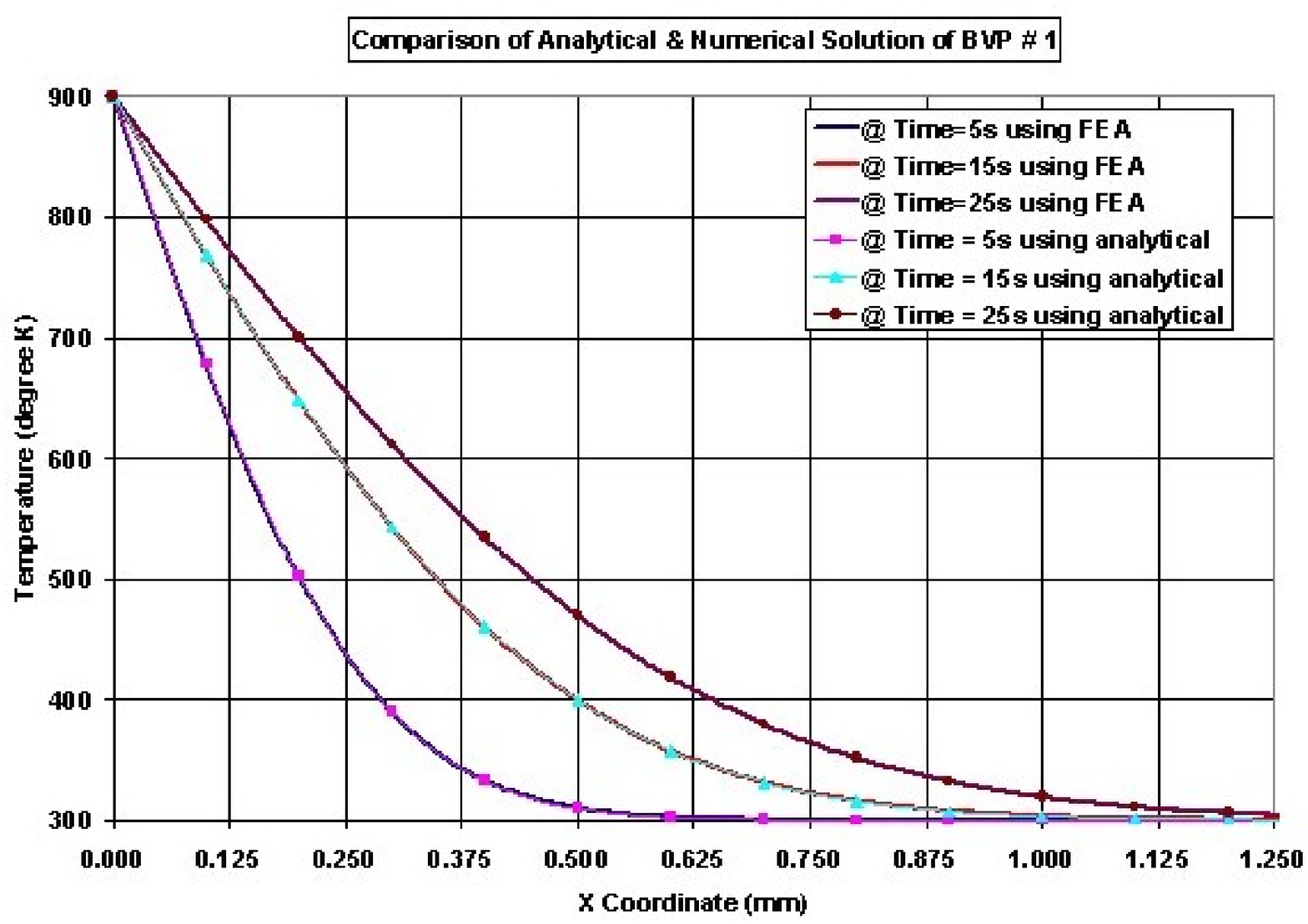}
  \caption{Comparison of analytical and numerical results for BVP-1}
  \label{fig:figure2}
\end{figure}

%\newpage
%

\begin{figure}[h] % float placement: (h)ere, page (t)op, page (b)ottom, other (p)age
  \centering
  % file name: D:/tahir-d/research/current-work/arif-bvps/first-paper/write-up/to-submit/arxiv/figure3.jpg
  \includegraphics[width=5.67in,height=3.95in,keepaspectratio]{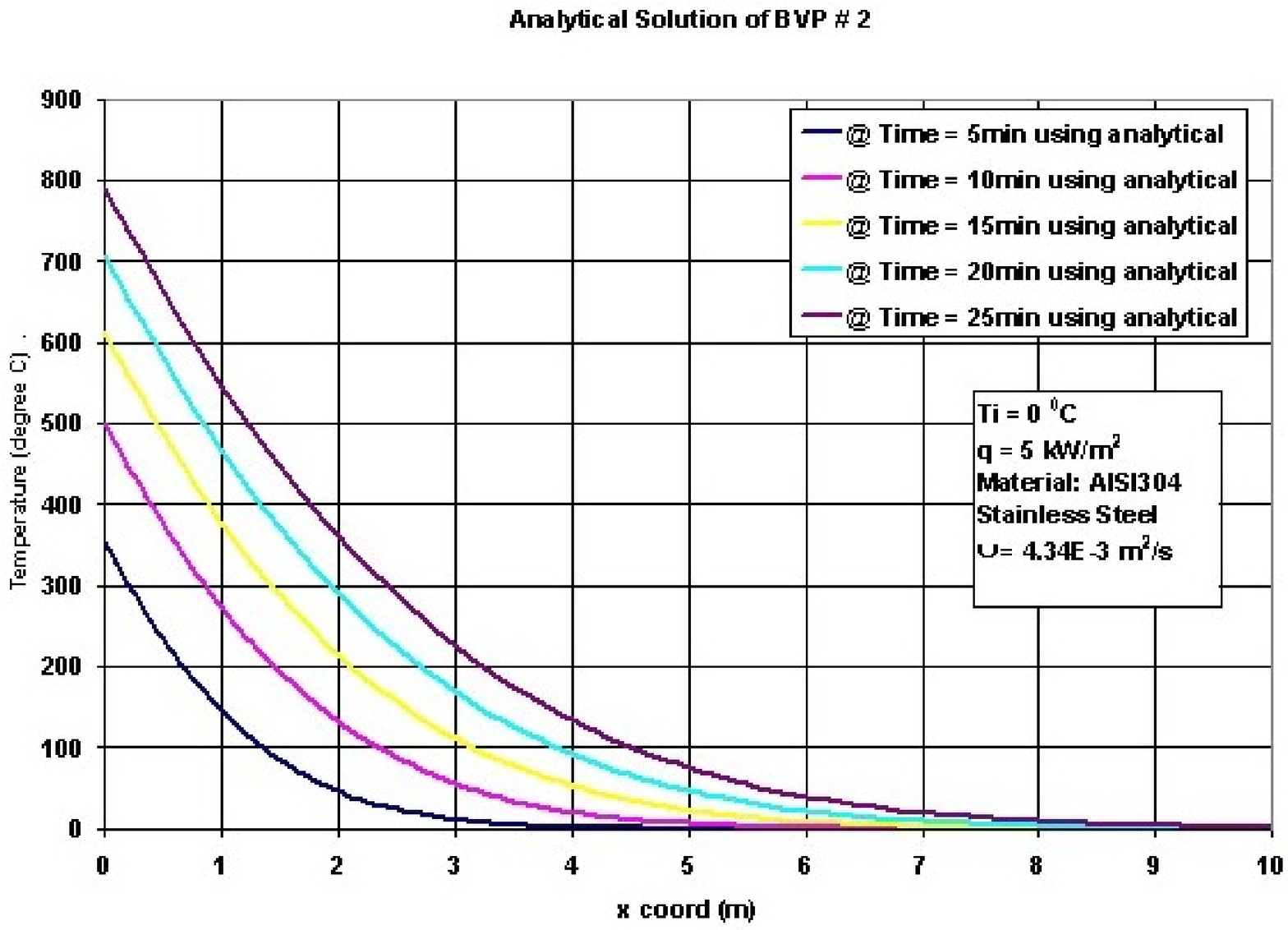}
  \caption{Temperature distribution at different times using analytical solution for BVP-2}
  \label{fig:figure3}
\end{figure}

\begin{figure}[h] % float placement: (h)ere, page (t)op, page (b)ottom, other (p)age
  \centering
  % file name: D:/tahir-d/research/current-work/arif-bvps/first-paper/write-up/to-submit/arxiv/figure4.jpg
  \includegraphics[width=5.67in,height=3.96in,keepaspectratio]{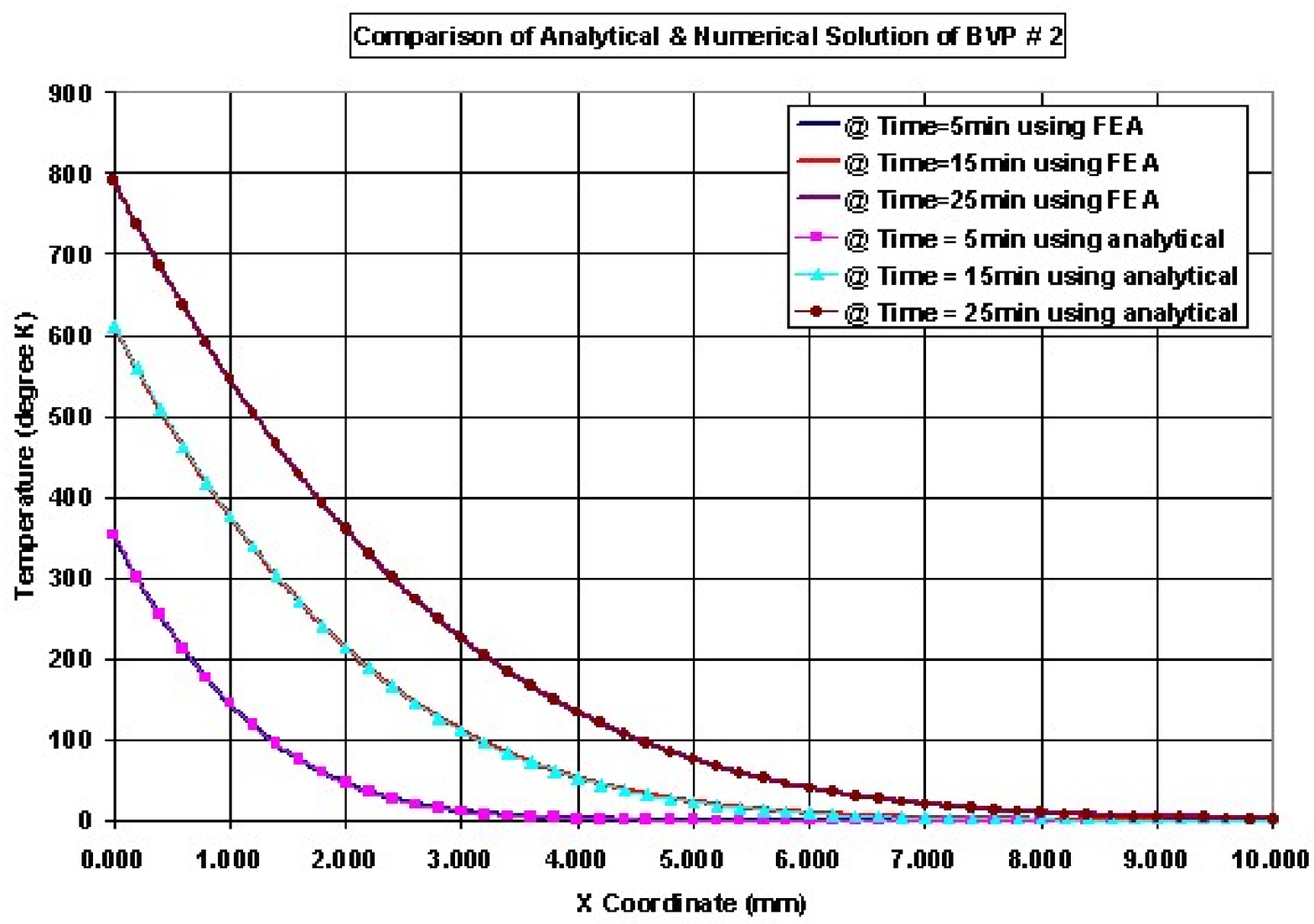}
  \caption{Comparison of analytical and numerical results for BVP-2}
  \label{fig:figure4}
\end{figure}

%\newpage

%%%%%%%%%%%%%%%%%%%%%
%%
%%
%%%%%%%%%%%%%%%%%%%%%

\ed